\documentclass{gtpart}     % Basic GT/GTM/AGT style

\title{Ideal boundary of 7-systolic complexes and groups}

\author{Damian Osajda}
\givenname{Damian}
\surname{Osajda}
\address{Instytut Matematyczny, Uniwersytet Wroc{\l}awski \newline
pl. Grunwaldzki 2/4, 50-384 Wroc{\l}aw, Poland \newline
and \newline
Institut de Math\'ematiques de Jussieu, Universit\'e Paris 6 \newline
Case 247, 4 Place Jussieu, 75252 Paris Cedex 05, France}
\email{dosaj@math.uni.wroc.pl}


\newtheorem*{wtw}{Theorem}
\newtheorem{lem}{Lemma}[section]
\newtheorem{tw}[lem]{Theorem}
\newtheorem{cor}[lem]{Corollary}
\newtheorem{prop}[lem]{Proposition}
\newtheorem*{rem}{Remark}
\newtheorem*{rems}{Remarks}
\theoremstyle{definition}
\newtheorem{de}[lem]{Definition}

\newcommand {\mr}{\mathrm}
\newcommand {\lk}{\left\{}
\newcommand {\rik}{\right\}}

\begin{document}

\begin{abstract}    % type your abstract below
We prove that ideal boundary of a $7$-systolic group is
strongly hereditarily aspherical. For some class of $7$-systolic
groups we show their boundaries are connected and without
local cut points, thus getting some results concerning splittings
of those groups.
\end{abstract}

\maketitle

%%%%%%%%%%%%%%%%%%%%   Start of main body of article

\section{Introduction}
The notion of $k$-systolic ($k\geq 6$ being a natural number)
complexes was introduced by T. Januszkiewicz and
J. \'Swi{\c a}tkowski
\cite{JS1} and, independently, by F. Haglund \cite{H} as
combinatorial analogue of nonpositively curved spaces.
Those complexes are simply connected simplicial complexes satisfying some
local combinatorial conditions. Roughly speaking there is a lower bound
for the length of ``essential" closed paths in a one-skeleton of every
link.

A group acting geometrically by automorphisms on a $k$-systolic
complex is called a $k$-systolic group. Examples of such torsion free groups
of arbitrary large cohomological dimension are constructed by Januszkiewicz and
\'Swi{\c a}tkowski \cite{JS1},
for every $k\geq 6$. Those examples are fundamental groups of some
simplices of groups. In the same paper it is proved that
$7$-systolic groups are Gromov-hyperbolic.

In this paper we study $7$-systolic complexes and groups and in particular
their ideal boundaries. Our main result is the following.

\begin{wtw}[Theorem \ref{tw1} in Section \ref{gro}]
\label{twA}
%\marginpar{[twA]}
The ideal boundary of a $7$-sys\-to\-lic
group is a strongly hereditarily aspherical compactum.
\end{wtw}

The notion of the strong hereditary asphericity (see Section
\ref{sha} for the precise definition) was introduced by R. J.
Daverman \cite{Da}. Roughly speaking a space is hereditarily
aspherical if each of its closed subsets is aspherical. The
significance of this notion follows from the fact that a cell-like
map defined on a strongly hereditarily aspherical compactum does
not raise dimension.

Theorem \ref{twA} shows that $7$-systolic groups are quite
different from many classical hyperbolic groups. It gives also new examples
of topological spaces that can occur as boundaries of hyperbolic groups.
The question about the spaces being such boundaries is well understood
only in dimensions $0$ and $1$---compare Kapovich and Kleiner \cite{KaKl}.
It is still not known
which topological spaces can be higher dimensional boundaries of
hyperbolic groups and only few homeomorphism types of such spaces are
known---see Benakli and Kapovich \cite{BK}
and discussion in p. 1) of Remarks in Section \ref{gro}.
Moreover, we show (see below) that for certain classes of complexes (and groups)
their ideal boundaries are ``simple'' in a sense---they are connected
and have no local cut points.
%From the resulta of Section \ref{ssc}
%it follows
%that all results concerning $7$-systolic complexes and groups
%apply also to simple convex cell complexes with $7$-large links.
%Those include e.g. $CAT(0)$ cube complexes
%with $7$-large links.

In order to prove the main theorem we define (in Section \ref{seven})
an inverse system of
combinatorial spheres in a $7$-systolic complex and projections onto
them, whose inverse limit
is the ideal boundary of the complex.
It should be noticed that even in the more general systolic (which means
$6$-systolic) case some inverse system like that can be defined
(cf. Januszkiewicz and \'Swi{\c a}tkowski \cite[Section 8]{JS1}).
However its properties
do not allow to consider it as a right tool to define
a reasonable boundary of a systolic complex or group. In particular
the inverse limit of this standard systolic inverse system
in a $7$-systolic case is not the Gromov boundary. On the other
hand, our construction is not valid in general systolic case,
although some results are probably true for hyperbolic
systolic (not necessarily $7$-systolic) groups.

In Section \ref{split} we study further properties of boundaries
of some $7$-systolic complexes. In particular we prove the following
theorem, which is a special case of Theorem \ref{tw3} proved
in that section.

\begin{wtw}[Corollary \ref{c9} in Section \ref{split}]
\label{twB}
%\marginpar{[twB]}
Let $X$ be a locally finite $7$-systolic normal pseudomanifold
of dimension at least $3$. Then its ideal boundary
$\partial X$ is connected and has no local cut points.
\end{wtw}

Via the results of Stallings \cite{S} and the ones of Bowditch
 \cite{B} the latter theorem implies the following.

\begin{wtw}[Corollary \ref{c10} in Section \ref{split}]
\label{twC}
%\marginpar{[twC]}
A group acting geometrically by automorphisms on a
locally finite $7$-systolic normal pseudomanifold
of dimension at least $3$ does not split essentially,
as an amalgamated product or as an
$HNN$-extension, over a finite nor two-ended group.
\end{wtw}

Groups acting geometrically on such pseudomanifolds of
arbitrary large dimension were constructed in Januszkiewicz and
\'Swi{\c a}tkowski \cite{JS1} and
are the only $7$-systolic groups of cohomological dimension
greater than $2$ known to us at the moment.

\subsection{Acknowledgements}
I would like to thank Tadeusz Januszkiewicz and Jacek \'Swi{\c a}tkowski for
introducing me to the subject and suggesting the problem. I am also
grateful to them and to Craig Guilbault,
Krzysztof Omiljanowski and Carrie Schermetzler for helpful
conversations.

Author was a Marie Curie
Intra-European fellow, contract MEIF CT 05011050 and was
partially supported by the Polish Scientific Research Committee
(KBN) grant 2 P03A 017 25

\section{Preliminaries}
\label{pre}
%\marginpar{[pre]}

\subsection{Simplicial complexes}
\label{sc}
%\marginpar{[sc]}
In this section we recall some definitions and fix the notation.

Let $X$ be a simplicial complex. We denote by $X'$ the first
barycentric subdivision of $X$. For a natural number $k$, we
denote by $X^{(k)}$ the {\it $k$-skeleton} of $X$, i.e. the union
of all simplices of $X$, of dimension at most $k$. For a given
subset $C=\lk v_1,v_2,...,v_l\rik$ of $X^{(0)}$ we denote by
$<v_1,v_2,...,v_l>$ the minimal simplex in $X$ containing
$C$---the simplex {\it spanned} by $C$. We denote by $X_{\sigma}$
the link of a given simplex $\sigma$ in $X$. A simplicial complex
$X$ is {\it flag} if every set $B$ of pairwise connected (by
edges) vertices of $X$ spans a simplex in $X$.

Recall that a subcomplex $Y$ of $X$ is {\it
full} if every set $B$ of vertices of $Y$ spanning a simplex of
$X$ spans a simplex in $Y$.
We denote
by $\sigma \ast
\rho$ the join of simplices $\sigma$ and $\rho$.

If it is not stated otherwise a simplicial complex $X$
is equipped with a path metric $d_X$ for which every $k$-simplex
of $X$ is isometric to the regular Euclidean $k$-simplex.

A simplicial complex $X$ is called a {\it chamber complex of
dimension $n$} if it is the union of $n$-simplices (which are
called {\it chambers} of $X$) and for every $(n-1)$-dimensional
face of $X$ there exist at least two chambers containing
that face. It is
easy to see that links in a chamber complex are themselves chamber
complexes. A {\it gallery} in a chamber complex is a finite
sequence of maximal simplices such that two consecutive simplices
share a common face of codimension $1$. A chamber complex is said
to be {\it gallery connected} if any two chambers can be connected
by a gallery. Chamber complex is {\it normal} if it is gallery
connected and all its links of dimension above $0$ are gallery
connected. A chamber complex is a {\it pseudomanifold} if every
codimension one face belongs to exactly two chambers.

\subsection{Systolic complexes and groups}
\label{scg}
%\marginpar{[scg]}

We follow here Januszkiewicz and \'Swi{\c a}tkowski
\cite{JS1}, \cite{JS2}.
For a given natural number $k\geq 4$, a simplicial complex $X$ is
{\it k-large} if it is flag and every cycle $\gamma$ in $X$ (i.e.
a subcomplex homeomorphic to the circle) of length $4\leq
|\gamma| <k$ has a diagonal (i.e. an edge connecting two
nonconsecutive vertices in $\gamma$). Here $|\gamma|$ denotes the
number of edges of $\gamma$.

A simplicial complex $X$ is {\it locally k-large} if for every
simplex $\sigma\neq \emptyset$ of $X$ its link $X_{\sigma}$ in $X$
is $k$-large.

$X$ is {\it k-systolic} if it is locally $k$-large, connected and
simply connected.

Because $k=6$ is of special importance in that theory,
$6$-systolic complexes are called {\it systolic}.

A group acting geometrically (i.e. properly discontinuously and
cocompactly) by simplicial automorphisms on a $k$-systolic (resp.
systolic) complex is called {\it k-systolic} (resp. {\it
systolic}). Free groups and fundamental groups of surfaces
are systolic groups. In Januszkiewicz and \'Swi{\c a}tkowski \cite{JS1},
for arbitrary $k$ and $n$,
a torsion free $k$-systolic group of cohomological dimension $n$
is constructed. Those groups are the fundamental groups of some
simplices of groups.

In the rest of this subsection we list some results concerning
systolic complexes and groups. We begin with the easy facts
whose proofs can be found in \cite{JS1}.

\begin{prop}
\label{p0.1}
%\marginpar{[p0.1]}
\begin{enumerate}
\item
If $k>m$ and $X$ is $k$-large then $X$ is also $m$-large.
\item
A full subcomplex in a (locally) $k$-large complex is (locally)
$k$-large.
\item
Links of a $k$-large complex are $k$-large.
\item
There is no $k$-large
triangulation of the $2$-sphere for $k\geq 6$. Hence no
triangulation of a manifold of dimension above $2$ is $6$-large
since $2$-spheres occur as links of some simplices in that case.
\end{enumerate}
\end{prop}

The following property of $7$-systolic complexes is crucial
for this paper.

\begin{tw} \cite[Theorem 2.1]{JS1}
\label{tw0.1}
%\marginpar{[tw0.1]}
Let $X$ be a $7$-systolic complex. Then the $1$-skeleton
of $X$ with its standard geodesic metric is
$\delta$-hyperbolic with $\delta=2\frac{1}{2}$.
\end{tw}

Thus $7$-systolic groups are word-hyperbolic.

In  \cite[Sections 3 and 7]{JS1} the notion of a
{\it convex} subcomplex of a systolic complex is
introduced. A simplex is a convex subcomplex.

%In Section 7 of \cite{JS1} it is also showed that
For a simplicial complex $X$ and its subcomplex $Q$ we can define
a {\it closed combinatorial ball of radius $i$ around $Q$ in $X$},
$B_i(Q,X)$, inductively: $B_0(Q,X)=Q$ and $B_i(Q,X)=\bigcup \lk
\tau : \tau \cap B_{i-1}(Q,X)\neq \emptyset \rik$, for every
positive natural number $i$.

By $S_i(Q,X)$ we denote the subcomplex of
$B_i(Q,X)$ spanned by the vertices at combinatorial distance
$i$ from $Q$, i.e. not belonging to $B_{i-1}(Q,X)$ (for $i>0$).
By ${\stackrel{\circ}{B_i}(Q,X)}$ we denote the {\it interior} of the closed
combinatorial $i$-ball around $\sigma$ in $X$, i.e.
${\stackrel{\circ}{B_i}(Q,X)}=B_i(Q,X)\setminus
S_i(Q,X)$.

For the rest of this section let $X$ denote a systolic
complex and $Q$ its convex subcomplex.

\begin{lem}\cite[Lemmas 7.5 and 7.6]{JS1}
\label{l0.2}
%\marginpar{[l0.2]}
The sphere $S_i(Q,X)$ and the ball $B_i(Q,X)$ are
full subcomplexes of $X$ and they are $6$-large.
\end{lem}

\begin{lem}\cite[Sect. 7]{JS1}
\label{l0.3}
%\marginpar{[l0.3]}
If $i>0$, then for every simplex $\tau
\in
S_i(Q,X)$, $\rho = \partial
B_{i-1}(Q,X)\cap X_{\tau}$ is a single simplex, $X_{\tau}
\cap B_{i}(Q,X)= B_1(\rho, X_{\tau})$
and $X_{\tau}
\cap S_{i}(Q,X)= S_1(\rho, X_{\tau})$.
\end{lem}

In the rest of the paper we call the simplex $\rho$, as in the
above lemma the
{\it projection of $\tau$ on $S_{i-1}(Q,X)$}.

The universal cover of a connected locally $6$-large simplicial
complex is systolic and the folowing fact holds.

\begin{tw}\cite[Theorem 4.1]{JS1}
\label{tw0.2}
%\marginpar{[tw0.2]}
The universal cover of a finite dimensional
connec\-ted locally $6$-large simplicial
complex is contractible.
\end{tw}

The proof of this theorem uses the projections onto
closed combinatorial balls (compare  \cite[Section 8]{JS1}).
Restrictions of those projections to spheres
$$
\pi_{B_i(Q,X)}|_{S_{i+1}(Q,X)}\colon S_{i+1}(Q,X)\to S_{i}(Q,X)
$$
have some properties that do not allow us to use them in
order to define a reasonable boundary of a systolic complex.
Thus in Section \ref{seven} we define, only for $7$-systolic
complexes, other maps between spheres.

The following lemma follows easily from the facts above and from
\cite[Corollary 1.5]{JS1}.

\begin{lem}
\label{l0.15}
%\marginpar{[l0.15]}
Let $k\geq 6$, let $Y$ be a $k$-large simplicial
complex and let $\sigma$ be a simplex of $X$.
If $p:X\to Y$ is the universal cover of $Y$ and
$m<\frac{k-1}{2}$ then for $i=0,1,2,...,m$ the map
$p|_{B_i(\sigma,X)}:B_i(\sigma,X)\to B_i(p(\sigma),Y)$ is an
isomorphism.
\end{lem}

We recall two results concerning systolic chamber complexes.

\begin{lem}\cite[Lemma 4.1]{O}
\label{l0.16}
%\marginpar{[l0.16]}
Let $X$ be a systolic chamber complex of dimension $n\geq 1$
and $\tau$ its simplex.
Then $S_k(\tau,X)$ is an $(n-1)$-dimensional chamber
complex, for every $k\geq 1$.
\end{lem}

\begin{lem}\cite[Lemma 4.2]{O}
\label{l0.17}
%\marginpar{[l0.17]}
Let $X$ be a systolic chamber complex of dimension $n\geq 1$
and $\tau$ its simplex.
Let $ \sigma$ be an $(n-1)$-dimensional simplex of
$S_k(\tau,X)$. Then there exists a vertex $v$ at a distance $k+1$ from
$\tau$ such that $v \ast \sigma$ is a simplex of $X$.
\end{lem}

\subsection{Strongly hereditarily aspherical compacta}
\label{sha}
%\marginpar{[sha]}
The notion of strongly hereditarily
aspherical compacta was introduced by R. J. Daverman
\cite{Da}. A compact metric space $Z$ is {\it strongly
hereditarily aspherical} if it can be embedded in the Hilbert cube
$Q$ in such a way that for each $\epsilon >0$ there exists an
$\epsilon$-covering ${\cal U}$ of $Z$ by open subsets of $Q$,
where the union of any subcollection of elements of ${\cal U}$ is
aspherical.

To show, in Section \ref{gro}, that the boundaries of $7$-systolic
groups are strongly hereditarily aspherical, we will use the
following result.

\begin{prop}\cite[Proposition 1]{Da}
\label{p1.1}
%\marginpar{[p1.1]}
Suppose $\lk L_i,\mu_i \rik$ is an
inverse sequence of finite polyhedra and PL maps, and suppose each
$L_i$ is endowed with a fixed triangulation $T_i$ such that
\begin{enumerate}
\item
$\mu_i^{-1}({\mr {each \; subcomplex \; of \;}}L_i)$ is
aspherical, and
\item
there exists a sequence $(a_k)_{k=1}^{\infty}$ of positive numbers,
such that ${\mr {lim}}_{k\rightarrow \infty}a_k=0$, and ${\mr
{diam}}(\mu_{i-k}\circ ... \circ \mu_{i-1}(\sigma))<a_k$, for
every simplex $\sigma \in T_i$.
\end{enumerate}
Then the inverse limit $Z={\mr {inv\; lim}}\lk L_i,\mu_i \rik$ is
strongly hereditarily aspherical.
\end{prop}

It should be noticed that, by Daverman and Dranishnikov
\cite[Theorem 2.10]{DD2}, every
strongly hereditarily aspherical compactum can be expressed as
an inverse limit like the one above.

\section{$7$-systolic complexes}
\label{seven}
%\marginpar{[seven]}
In this section we study some properties of
$7$-systolic complexes. In particular we define and examine other
(then in the general systolic case) contractions on spheres. This
is crucial for Section \ref{gro}. Then we study properties of
some $7$-systolic chamber complexes. Those results are important
for Section \ref{split}.

Let $X$ be a $7$-systolic complex of dimension $n<\infty$. Let $Q$
be its convex subcomplex (see
\cite[Sections 3 and 7]{JS1} and
compare Section \ref{scg} above). For a natural number $k$, we
denote by $S_k$ the combinatorial sphere $S_k(Q,X)$ (compare
Section \ref{scg} above) and we denote by $B_k$ the closed ball
$B_k(Q,X)$.

Define a map $\pi_{Q,k}\colon S_k^{(0)}\to (S_{k-1}')^{(0)}$ by
putting $\pi_{Q,k}(w)=b_{\tau}$, for every vertex $w$ of $S_k$,
where $\tau =X_w\cap B_{k-1}$ is the projection of $w$ on
$S_{k-1}$ and $b_{\tau}\in (S_{k-1}')^{(0)}$ is the barycenter of
$\tau$.

\begin{lem}
\label{l1}
%\marginpar{[l1]}
Let $v_1$ and $v_2$ be two
vertices in $S_k$ belonging to the same simplex. Then
$\pi_{Q,k}(v_1)$ and $\pi_{Q,k}(v_2)$ belong to the same simplex
of $S_{k-1}'$.
\end{lem}
\begin{proof} Let $\tau=X_{<v_1,v_2>}\cap B_{k-1}$ and $\tau_i=X_{v_i}\cap
B_{k-1}$, for $i=1,2$. Then $\tau\subset \tau_1 \cap \tau_2$. It
is enough to show that $\tau_1\subset \tau_2$ or $\tau_2\subset
\tau_1$.

We will show this arguing by contradiction. Suppose it is not
true, i.e. there exist vertices $w_i$ such that $w_i\in \tau_i
\setminus\tau_j$ for $\lk i,j\rik =\lk 1,2\rik$. Let for $i=1,2$,
$t_i$ be a vertex belonging to $X_{w_i}\cap B_{k-2}\cap
X_{\tau}$.

Let us examine the closed path $(v_1,w_1,t_1,t_2,w_2,v_2,v_1)$ in
$X^{(1)}$.

There are no diagonals of the form $<v_i,t_j>$ since the distance
between $v_i$ and $t_j$ is at least $2$, $i,j=1,2$.

There are no diagonals of the form $<v_i,w_j>$, $i\neq j$. Since
if, e.g. $<v_1,w_2>$ is an edge in $X$ then $w_2\in X_{v_1}\cap
B_{k-1}=\tau_1$.

There is no diagonal $<w_1,w_2>$. If it exists, then the path
$(v_1,w_1,w_2,v_2,v_1)$ is a closed simple path without diagonals
of length $4$ which contradicts $7$-largeness of $X$.
Similarily there are no diagonals $<w_1,t_2>$ and $<w_2,t_1>$.

Hence the path $(v_1,w_1,t_1,t_2,w_2,v_2,v_1)$ is a closed path of
length at most six without diagonals. This contradicts
$7$-largeness of $X$.
\end{proof}

Using the above lemma we can extend $\pi_{Q,k}$ simplicially.

\begin{de}
\label{d1}
%\marginpar{[d1]}
Define, for a natural number $k$, a contiunous map between
combinatorial spheres
$$
\pi_{Q,k}\colon S_k(Q,X)\to
(S_{k-1}(Q,X))',
$$
given by the simplicial extension
of the map
$$\pi_{Q,k}\colon S_k(Q,X)^{(0)}\to
(S_{k-1}(Q,X)')^{(0)}.
$$
\end{de}

\begin{lem}
\label{l1a}
%\marginpar{[l1a]}
There exists a constant $C<1$, depending only on
$n={\mr {dim}}(X)$ such that for every $k,l\in \lk1,2,3,...\rik$ with
$l<k$ and for every two points $x,y\in S_k$ one has ${\mr
d}_{S_{k-l-1}}(\pi_{Q,k-l}\circ ... \circ \pi_{Q,k-1}\circ
\pi_{Q,k}(x),\pi_{Q,k-l}\circ ... \circ \pi_{Q,k-1}\circ
\pi_{Q,k}(y))\leq C^{l+1} {\mr d}_{S_k}(x,y)$.
\end{lem}
\begin{proof}
Let $D$ be the distance from a vertex to an opposite
codimension one face in the regular $(n-1)$-simplex. Let $E$
be the diameter of a maximal simplex in
the barycentric subdivision of the regular $(n-1)$-simplex.
Then for $C=\frac{E}{D}<1$ the lemma holds.
\end{proof}

\begin{lem}
\label{l2}
%\marginpar{[l2]}
For every subcomplex $L$ of $S_{k-1}$ the subcomplex
$\pi_{Q,k}^{-1}(L)$ of $S_k$ is aspherical.
\end{lem}
\begin{proof} We will show that $\pi_{Q,k}^{-1}(L)$ is a full subcomplex of
$S_k$. Hence, by Proposition \ref{p0.1}, as a full subcomplex of
$6$-large complex it is $6$-large and thus aspherical,
by Theorem \ref{tw0.2}.

Let vertices $v_1,v_2,...,v_l\in (\pi_{Q,k}^{-1}(L))^{(0)}$ span a
simplex in $S_k$. Then, by Lemma \ref{l1},
$\pi_{Q,k}(v_1),...,\pi_{Q,k}(v_l)$ are vertices of a simplex of
$S_{k-1}'$ and they correspond to a chain of simplices
$\tau_1,...,\tau_l$ of $S_{k-1}$. One of them, say $\tau_1$, is the
highest dimensional simplex among them and hence it contains
all the simplices $\tau_2,...,\tau_l$. This means that
$\pi_{Q,k}^{-1}(\tau_1)\subset \pi_{Q,k}^{-1}(L)$ contains all
points $v_1,...,v_l$ and hence simplex spanned by them.
\end{proof}

In the rest of this section we study some $7$-large chamber
complexes.

For a $7$-large chamber complex $X$ we define, for a vertex $v$ of $X$, a
condition $R(v,X)$ that will be crucial for Section \ref{split}
(compare the condition $R(v,X)$ defined in \cite[Section 4]{O}):
%\\
%\hspace{2cm}
$$
R(v,X)\; \; {\mr{iff}}\; \; \; (\forall \; \sigma \in X_v \;(X_v \setminus
{\stackrel{\circ}{B_2}}(\sigma,X_v) {\rm \; and \;} X_v \setminus
{\stackrel{\circ}{B_3}}(\sigma,X_v){\rm \; are \; connected}))
$$

The next lemma is an analogue of \cite[Lemma 4.7]{O}, for
$7$-large complexes.

\begin{lem}
\label{l0.5}
%\marginpar{[l0.5]}
Let $X$ be a $7$-large chamber complex such that the
link $X_{\kappa}$ is connected, for every simplex $\kappa$ of $X$
of codimension greater than one, and $X_{\sigma}\setminus
\stackrel{\circ}{B_i}(\rho,X_{\sigma})$, $i=2,3$ is connected for every
codimension two simplex $\sigma$ of $X$ and every simplex $\rho$
of its link $X_{\sigma}$. Then for every vertex $v$ of $X$
condition $R(v,X)$ holds.
%Then$X_{v}\setminus B_i(\omega,X_v)$, $i=2,3$ is connected for
%every vertex $v$ of $X$and every simplex $\omega$ of its link
%$X_v$
\end{lem}
\begin{proof} We will proceed by induction on $n={\rm dim}(X)$.

For $n=2$ the assertion is clear since codimension two simplexes
are just vertices.

Assume we proved the lemma for $n\leq k$. Let ${\rm dim}(X)=k+1$.
Take a vertex $v$ of $X$ and consider its link $X_v$. It has
dimension $k$. Moreover for every codimension $l$ simplex $\sigma$
of $X_v$ the simplex $\sigma \ast v$ is of codimension $l$ in $X$
and $X_{\sigma \ast v}=(X_v)_{\sigma}$. Thus $X_v$ satisfies
hypotheses of the lemma. Hence by the induction assumptions, for
every vertex $w$ of $X_v$ condition $R(w,X_v)$ holds.

$X_v$ is $7$-large as a full subcomplex of
$X$ (Proposition \ref{p0.1}).
Thus the universal cover $\widetilde{X_v}$ of $X_v$ is $7$-systolic.
Let $p: \widetilde{X_v}\rightarrow
X_v$ be a covering map.

Take a simplex $\omega$ of $X_v$.
Let $\kappa$ be a simplex of
$\widetilde{X_v}$
such that $p(\kappa)=\omega$.
Since
$S_0(\kappa,\widetilde{X_v})$ is connected we have,
by \cite[Lemma 4.5]{O}
(or by Corollary \ref{c8} below), that
$S_1(\kappa,\widetilde{X_v})$, $S_2(\kappa,\widetilde{X_v})$
and $S_3(\kappa,\widetilde{X_v})$ are connected.
By Lemma \ref{l0.15},
$p'=p|_{B_2(\kappa,\widetilde{X_v})}:
{B_2(\kappa,\widetilde{X_v})}\to {B_2(\omega,X_v)}$
is an isomorphism.
%$p({B_2(\kappa,\widetilde{X_v})})=
%{B_2(\omega,X_v)}$.
Observe that $p({B_3(\kappa,\widetilde{X_v})})\subset
{B_3(\omega,X_v)}$. We want to show that
$p({B_3(\kappa,\widetilde{X_v})})=
{B_3(\omega,X_v)}$. Let $z$ be a vertex in
${B_3(\omega,X_v)}\setminus {B_2(\omega,X_v)}$ and let
$u\in {B_2(\omega,X_v)}$ be a vertex connected by an edge
with $z$. Then, by Lemma \ref{l0.15},
$p''=p|_{B_1((p')^{-1}(u),\widetilde{X_v})}:
{B_1((p')^{-1}(u),\widetilde{X_v})}\to {B_1(u,X_v)}$
is an isomorphism and $(p'')^{-1}(z)\in
{B_3(\kappa,\widetilde{X_v})}$. Hence
$z\in p({B_3(\kappa,\widetilde{X_v})})$ and
$p({B_3(\kappa,\widetilde{X_v})})=
{B_3(\omega,X_v)}$.
Now we claim that $S_3(\omega,X_v)=
p(S_3(\kappa,\widetilde{X_v}))$ and hence is connected.
Observe that $S_3(\omega,X_v)\subset
p(S_3(\kappa,\widetilde{X_v}))$. Suppose
$S_3(\omega,X_v)\neq
p(S_3(\kappa,\widetilde{X_v}))$. Let $w_1\in
S_3(\kappa,\widetilde{X_v})$ be a vertex such that
$p(w_1)\in B_2(\omega,{X_v})$. The vertex
$w_2=(p')^{-1}(p(w_1))$ belongs to
$B_2(\kappa,\widetilde{X_v})$ and $p(w_1)=p(w_2)$.
But then
$d_{\widetilde {X_v}^{(1)}}(w_1,w_2)<7$ and we can find
homotopically non-trivial closed path of length less than $7$ in
$X_v$. This contradicts $7$-largeness of $X_v$, by
\cite[Corollary 1.5]{JS1}. Thus we have shown that
$S_3(\omega,X_v)=
p(S_3(\kappa,\widetilde{X_v}))$ is connected.

Take two vertices $t$ and $s$ of $X_v\setminus
\stackrel{\circ}{B_2}(\omega,X_v)$ (or of $X_v\setminus
\stackrel{\circ}{B_3}(\omega,X_v)$).
Since, by assumptions, $X_v$ is connected there exists a path in
$(X_v)^{(1)}$ joining them. If this path misses
$\stackrel{\circ}{B_2}(\omega,X_v)$ (respectively $\stackrel{\circ}
{B_3}(\omega,X_v)$) it
joins these vertices in $X_v\setminus \stackrel{\circ}
{B_2}(\omega,X_v)$
(respectively in $X_v\setminus \stackrel{\circ}
{B_3}(\omega,X_v)$). If not we
can replace it, by connectedness of $S_2(\omega,X_v)$
(respectively $S_3(\omega,X_v)$), by a path intersecting
$S_2(\omega,X_v)$ (respectively $S_3(\omega,X_v)$) and also lying
in $X_v\setminus \stackrel{\circ}
{B_2}(\omega,X_v)$ (respectively in
$X_v\setminus \stackrel{\circ}{B_3}(\omega,X_v)$).
Hence we get the
conclusion. \end{proof}

\begin{cor}
%[c0.6]
\label{c0.6} Let $X$ be a normal $7$-systolic pseudomanifold. Then
condition $R(v,X)$ holds for every vertex $v$ of $X$.
%Then$X_{v}\setminus B_i(\omega,X_v)$, $i=2,3$ is connected for
%every vertex $v$ of $X$and every simplex $\omega$ of its link
%X_v$.
\end{cor}

\begin{proof} One dimensional link in a normal pseudomanifold is a circle.
Hence it satisfies assumptions of the preceding corollary. \end{proof}

\section{Gromov boundary}
\label{gro}
%\marginpar{[gro]}
The aim of this section is to prove that the ideal
boundary of a $7$-systolic group (such groups are word hyperbolic
by Theorem \ref{tw0.1}) is a strongly hereditarily
aspherical compactum. To prove this we first show that such a
boundary can be described as an inverse limit of combinatorial
spheres in the complex on which the group acts geometrically.

Throughout this section $X$ denotes a locally finite $7$-systolic
complex of dimension $n<\infty$. We fix a vertex $v$ of $X$. For a
natural number $k$, we denote by $S_k$ the combinatorial sphere
$S_k(v,X)$ and we denote by $B_k$ the closed ball $B_k(v,X)$. We
denote by $\pi_k$ the projection $\pi_{\lk v \rik,k}\colon S_k\to
S_{k-1}$ (see Section \ref{seven}).
%Let moreover $\delta X={\mr
%{inv \; lim}}(S_k, \pi_k)$.

\begin{lem}
\label{l3}
%\marginpar{[l3]}
$\delta X={\mr
{inv \; lim}}\lk  S_k, \pi_k\rik$ is homeomorphic to
$\partial X$--the Gromov boundary of $X$.
\end{lem}
\begin{proof} We use the set
of equivalence classes of geodesic rays in $X^{(1)}$ propagating
from a given vertex $v$, as a definition of the Gromov boundary
of $X$---for details see e.g.
Bridson and Haefliger \cite[Chapter III.3]{BH}.

Compactness of both $\delta X$ and $\partial X$ follows from the
fact that the balls in $X^{(0)}$ are finite.

First, we construct a bijection $F\colon \delta X \to \partial X$.
Let $x=(v,x_1,x_2,...)\in \delta X$. Note that $x_k\in S_k$ for
$k=1,2,...$. For arbitrary $k$, choose a maximal simplex
$\sigma_k$ of $S_k$ containing $x_k$. If we take a vertex $u$ of
$\sigma_k$ then $X_{\sigma_k}\cap B_{k-1}\subset X_u\cap B_{k-1}$
and hence there exists a vertex $w$ of $\sigma_{k-1}$ connected by
an edge with $u$. Hence for any $k$ we can construct a sequence
$(v=v_0^k,v_1^k,v_2^k,...,v_k^k)$ of vertices of $X$ such that
$v_i^k\in \sigma_i$ and $v_k^{i-1}$ is connected by an edge with
$v_k^i$ for $i=1,2,...,k$. Since a path $c_k=(v_0^k,...,v_k^k)$
in $X^{(1)}$ has length $k$ and joins $v$ and $v_k^k$ lying at a
distance $k$ it is a geodesic segment starting at $v$. Since 
balls in $X^{(0)}$ are finite, we can, by the diagonal argument,
extract from $(c_k)_{k=1}^{\infty}$ a geodesic ray
$c=(v,v_1,v_2,...)$, such that $v_k$ is a vertex of $\sigma_k$,
for every $k=1,2,3,...$. The equivalence class of $c$ within
$\partial X$ is by definition $F(x)$. Observe that it is
independent of choosing $c$ as above, since all of them lie at
distance at most one from the sequence $x$.

We show now that $F$ is injective. Let $x=(v,x_1,x_2,...)$ and
$y=(v,y_1,y_2,...)$ be two elements of $\delta X$ with
$F(x)=F(y)$. Let the geodesic rays $c=(v,v_1,v_2,...)$ and
$d=(v,w_1,w_2,...)$ in $X^{(1)}$ representing, respectively,
$F(x)$ and $F(y)$ be constructed as above. Then there exists a
constant $D>0$ such that for every $k=1,2,3,...$ we have
$d_{S_k}(v_k,w_k)\leq D$. Fix $k$. It is enough to show that
$d_{S_k}(x_k,y_k)\leq \epsilon$ for every $\epsilon > 0$. Choose
$\epsilon >0$. Take $l\in {\bf N}$ such that $l\geq {\mr
{log}}_C\frac{\epsilon}{D+2}$, where $C<1$ is the constant from the
proof of Lemma
\ref{l1a}. By construction $d_{S_{k+l}}(x_{k+l},y_{k+l})\leq D+2$
and thus by Lemma \ref{l1a}
$$
d_{S_k}(x_k,y_k)=
$$
$$
=d_{S_k}(\pi_{k+1}\circ...\circ
\pi_{k+l-1}\circ \pi_{k+l}(x_{k+l}),\pi_{k+1}\circ...\circ
\pi_{k+l-1}\circ \pi_{k+l}(y_{k+l}))\leq
$$
$$
\leq C^{l}d_{S_{k+l}}(x_{k+l},y_{k+l})\leq \epsilon.
$$

Now we show $F$ is onto. Take a geodesic ray $c=(v,v_1,v_2,...)$,
$v_i\in X^{(0)}$. Observe that $v_k\in S_k$ for $k=2,3,4,...$.
Consider a sequence $(\pi_2\circ \pi_3\circ...\circ
\pi_k(v_k))_{k=2}^{\infty}$ of points in $S_1$. By compactness of
spheres there is a subsequence $(v_{a^1(1)},v_{a^1(2)},...)$ of
the sequence $(v_2,v_3,...)$ such that $(\pi_2\circ
\pi_3\circ...\circ \pi_{a^1(k)}(v_{a^1(k)}))_{k=2}^{\infty}$
converges. Let $x_1\in S_1$ be the limit of this subsequence. Now
given a subsequence (of $(v,v_1,v_2,...)$)
$(v_{a^l(1)},v_{a^l(2)},...)$, $l>1$ we find a subsequence
$(v_{a^{l+1}(1)},v_{a^{l+1}(2)},...)$, $a^{l+1}(i)> l$ such that
$\pi_{l+2}\circ...\circ
\pi_{a^{l+1}(k)}(v_{a^{l+1}(k)}))_{k=l+1}^{\infty}$ tends to
$x_{l+1}\in S_{l+1}$. By construction $\pi_k(x_k)=x_{k-1}$, for $k>1$ and
$\pi_1(x_1)=v$. Hence $x=(v,x_1,x_2,...)\in \delta X$. Moreover
since $v_k$ and $v_{k+1}$ belong to a common simplex for every
$k=2,3,4...$ we get, by definition of $\pi_k$, that
$d(v_k,\pi_{k+1}(v_{k+1}))\leq E$, for $E<1$ being the constant
from the proof of Lemma \ref{l1a}. 
Then for every natural number $l$ we have
$d(v_k,\pi_{k+1}\circ \pi_{k+2} \circ...\circ
\pi_{k+l}(v_{k+l}))\leq \sum_{i=1}^{\infty}E^i<\infty$. Thus
$d(v_k,x_k)< \sum_{i=1}^{\infty}E^i<\infty$, that implies $c$
represents the equivalence class of $F(x)$.

Finally we argue $F$ is continuous and hence as a continuous
bijection defined on a compact space it is a homeomorphism.

Given $x=(v,x_1,x_2,...)\in \delta X$ and a sequence
$(x^i)_{i=1}^{\infty}\subset \delta X$, $x^i=(v,x_1^i,...)$
converging to $x$, fix geodesic rays $c=(v,v_1,v_2,...)$ and
$c_i=(v,v_1^i,v_2^i,...)$, $i=1,2,3,...$ representing,
respectively, $F(x)$ and $F(x^i)$, $i=1,2,3,...$, and constructed
as when we defined $F$. To prove $F$ is continuous at $x$ it is
enough to show that for every natural number $N$ there exists
$M>0$ such that for every natural number $i>M$ we have
$d_{S_N}(v_N,v_N^i)<3$. By definition of the topology of an
inverse limit there exists $M>0$ such that for every natural
number $i>M$ one has $d_{S_N}(x_N,x_N^i)<1$ and hence
$d_{S_N}(v_N,v_N^i)\leq d_{S_N}(v_N,x_N)+ d_{S_N}(x_N,x_N^i)
+d_{S_N}(x_N^i,v_N^i)<3$. \end{proof}

Now we state and prove the following main theorem.

\begin{tw}
\label{tw1}
%\marginpar{[tw1]}
The ideal boundary of a $7$-systolic
group is a strongly hereditarily aspherical compactum.
\end{tw}
\begin{proof}
A $7$-systolic group $G$ acts, by definition, geometrically on a locally finite
$7$-systolic complex $X$ of finite dimension. Then the ideal boundary
$\partial G$ of
$G$ is homeomorphic to $\partial X$.

We apply Proposition \ref{p1.1} to the inverse system
$\lk L_i,\mu_i \rik=\lk S_i,\pi_{i+1}\rik$. By Lemma \ref{l2} the condition
$1)$ of Proposition \ref{p1.1} is fulfilled, and by Lemma \ref{l1a} we get
condition $2)$ of the proposition. Hence $\partial G=\partial X=\delta X=
{\mr {inv\; lim}}\lk S_i,\pi_{i+1} \rik$ is a strongly hereditarily aspherical
compactum.
\end{proof}

\begin{rems}
{\bf 1)} A simple argument shows that every compact metrizable space
can be homeomorphic to the ideal boundary of some hyperbolic space
(even more---of some $CAT(-1)$ space).
The question of which topological
spaces can occur as boundaries of hyperbolic
groups (compare Benakli and Kapovich \cite[Chapter 17]{BK}) is more difficult.
It is answered somehow only in dimensions (of the boundary) $0$ and
$1$ (cf. Kapovich and Kleiner \cite{KaKl}). For higher dimensions the following homeomorphism
types of the boundaries of hyperbolic groups were known:
spheres,
Pontryagin surfaces $\Pi_p$ for $p$ being a prime number, two-dimensional
universal Menger compactum $\mu_2^5$ (compare \cite[Chapter 17]{BK}),
three-dimensional universal Menger compactum $\mu_3^7$
(cf. Dymara and Osajda \cite{DO}), Pontryagin spheres and three-dimensional trees
of manifolds (cf. Przytycki and \'Swi{\c a}tkowski \cite{PS}).

By Theorem \ref{tw0.1}, $7$-systolic groups are hyperbolic and, by
\cite[Corollary 19.3]{JS1}, for each natural number $n$, there
exists a hyperbolic group of cohomological dimension $n$. Hence,
by Theorem \ref{tw1} above, and by Bestvina and Mess
\cite[Corollary 1.4]{BesM},
strongly hereditarily aspherical compacta of all dimensions can
occur as boundaries of hyperbolic groups.

Moreover, in \cite{JS1} examples of $7$-systolic groups acting
on pseudomanifolds of arbitrary large dimension are constructed.
Thus, by Corollaries \ref{c9} and \ref{c10},
those group are, in a sense, indecomposable and their boundaries
are connected, locally connected and without local cut points (compare
Section \ref{split}).

{\bf 2)} Zawi\'slak \cite{Z} has shown that the
boundary of a $7$-systolic orientable normal pseudomanifold of
dimension $3$
is the Pontryagin sphere (cf. Jakobsche \cite{Jak}).
Such pseudomanifolds are constructed in
\cite{JS1}.

%\vspace{0.3cm}
%\epsfxsize=11cm
%\frame{\epsffile{ponswia.ps}}
%\begin{center} Approaching Pontryagin sphere
%\vspace{0.3cm}

The Pontryagin sphere
%\marginpar{referencje ?}
is the inverse limit of an inverse system
$\lk X_i,p_i \rik_{i=1}^{\infty}$ defined as follows. Let $X_1=S^2$
be a triangulated two-sphere. Let ${\cal T}$ be a given
triangulation of the two torus. Assume $X_i,p_j$ are defined for
$i\leq k$ and $j\leq k-1$. Let $X_k$ be a surface and ${\cal T}_k$ its
triangulation. $X_{k+1}$ is a connected sum of $X_k$ and a set
of disjoint
tori $T_{\sigma}$--one for every $2$-simplex $\sigma$ of
${\cal T}_k$--carrying the triangulation ${\cal T}$.
Every $T_{\sigma}$ is glued to $X_k$
by identifying $\partial \sigma$
and the boundary of some $2$-simplex $\sigma '$
of triangulation of $T_{\sigma}$.
Then $X_{k+1}$ carries an induced triangulation and we define
a triangulation ${\cal T}_{k+1}$ of $X_{k+1}$ as a subdivision
of this natural triangulation. The map $p_k\colon X_{k+1}\to X_k$
is defined by the conditions: $p_k(T_{\sigma}\setminus \sigma ')=
{\mr {int}\; }{\sigma}$ and $p_k|_{\partial \sigma}=
Id_{\partial \sigma}$ for every $2$-simplex $\sigma$ of
${\cal T}_k$.

{\bf 3)} For a polytopal complex $Y$ its {\it face complex} $\Phi(Y)$
is a simplicial complex defined as follows. The vertex set of $\Phi(Y)$
is the set of cells of $Y$ and the vertices of $\Phi(Y)$ span a simplex
if the cells of $Y$ corresponding to those vertices are contained in
a common cell of $Y$. It can be shown (compare \cite{JS3}) that if
$Y$ is a simply connected simple (i.e. all links are simplicial complexes)
polytopal complex with $7$-large links then $\Phi(Y)$ is $7$-systolic.
Thus the ideal boundary of such a complex $Y$ is strongly hereditarily
aspherical.
\end{rems}
{\bf Question.}

Is the ideal boundary of a hyperbolic systolic group strongly
hereditarily aspherical ?

\section{Splittings}
\label{split}
The aim of this section is to study further properties
of boundaries of $7$-systolic complexes in some special cases.
As a consequence we get results concerning splittings
of groups acting on such complexes.

Throughout this section $X$ denotes a locally finite $7$-systolic
chamber complex of dimension $n<\infty$. We fix a vertex $v$ of $X$. For a
natural number $k$, we denote by $S_k$ the combinatorial sphere
$S_k(v,X)$ and we denote by $B_k$ the closed ball $B_k(v,X)$. We
denote by $\pi_k$ the projection $\pi_{\lk v \rik,k}\colon S_k\to
S_{k-1}$ (see Section \ref{seven}).

\begin{lem}
\label{l4}
%\marginpar{[l4]}
Let $Y$ be a $7$-large
$n$-dimensional chamber complex, $\sigma$ one of its simplices
and
$\tau$ an $(n-1)$-simplex of $S_2(\sigma,Y)$.
Then there
exists a vertex $v\in Y\setminus B_2(\sigma,Y)$ such
that $v\ast \tau$ is a simplex of $Y$.
\end{lem}
\begin{proof}
By Lemma \ref{l0.17}, if
we consider the universal cover $p:\widetilde Y\rightarrow Y$ and
$\widetilde {\sigma}\in p^{-1}(\sigma)$, $\widetilde {\tau}\in
p^{-1}(\tau)\cap S_2(\widetilde {\sigma},\widetilde Y)$, then
there exists a vertex $\widetilde v$ of $\widetilde Y$ such that
$\widetilde
v\in \widetilde Y\setminus B_2(\widetilde {\sigma},\widetilde Y)$
and $\widetilde v \ast\widetilde \tau$ is a simplex of $\widetilde Y$.
Consider $v=p(\widetilde v)$. Clearly $v\ast \tau$ is a simplex
of $Y$.

Assume $v\in B_2(\sigma,Y)$. Then there exists a simplex
$\widetilde {\sigma_1}\in p^{-1}(\sigma)$ distinct from $\widetilde
{\sigma}$ such that $\widetilde v\in B_2(\widetilde
{\sigma_1},\widetilde Y)$. Since $\widetilde v \in B_3(\widetilde
{\sigma},\widetilde Y)$ we can then choose vertices $s\in \widetilde
{\sigma}$ and $t\in \widetilde {\sigma_1}$ with $p(s)=p(t)$ and a
path of length at most $6$ joining $s$ and $t$. But this
contradicts $7$-largeness of
$Y$. Thus $v\in Y\setminus B_2(\sigma,Y)$.
\end{proof}

\begin{lem}
\label{l6}
%\marginpar{[l6]}
The map $\pi_k\colon S_k \to S_{k-1}$ is onto.
\end{lem}
\begin{proof}
Let $z$ be a given point in $S_{k-1}$. We will show that
there exists a
point $w\in S_k$ satisfying $\pi_k(w)=z$.

{\it Case 1:} $z$ is a
barycenter of a simplex $\tau$ of $S_{k-1}$.

If ${\mr
{dim}}(\tau)=n-1$ then by Lemma \ref{l0.17} there
exists a vertex $w\in S_k$
such that $w\ast \tau$ is a simplex of $X$ and hence
$\pi_k(w)=z$.

Now, let ${\mr {dim}}(\tau)=m<n-1$. Since, by Lemma \ref{l0.16},
$S_{k-1}$ is a chamber complex of dimension $n-1$, there
exists an $(n-1)$-simplex
$\rho$ of $S_{k-1}$ containing $\tau$. Then, again by
Lemma \ref{l0.17},
there
exists a vertex $w'\in S_k$ spanning a simplex with $\rho$.
Clearly $w'\in S_k \cap S_2(\delta,X_{\tau})$,
where $\delta=X_{\tau}\cap S_{k-2}$. Since $X_{\tau}$ is a
$7$-large $(n-m-1)$-dimensional chamber complex and,
(again by Lemma \ref{l0.16})
$S_2(\delta,X_{\tau})$ is an $(n-m-2)$-dimensional chamber
complex (nonempty), we get, by Lemma \ref{l4}, that there exists
a vertex $w\in X_{\tau}\setminus B_2(\delta,X_{\tau})$. It follows
that
$\pi_k(w)=z$.

{\it Case 2:} $z$ belongs to an interior of an
$m$-simplex $\tau$ of $S_{k-1}'$.

Then $\tau=<a_0,a_1,...,a_m>$
where $a_i$ is the barycenter of an $i$-simplex $\tau_i$ of
$S_{k-1}$. By {\it Case 1} there exists a vertex $a_m'\in S_k$ such that
$\pi_k(a_m')=a_m$. Then $a_m'\in S_2(X_{\tau_{m-1}}\cap
{B_{k-2}},X_{\tau_{m-1}})$ and, using Lemma \ref{l4}
for $ X_{\tau_{m-1}}$, there exists a vertex $a_{m-1}'\in
X_{\tau_{m-1}}\setminus {B_2( X_{\tau_{m-1}}\cap
{B_{k-2}},X_{\tau_{m-1}})}$ connected to $a_m'$ by an
edge. Note that $\pi_k(a_{m-1}')=a_{m-1}$ and that
$<a_m',a_{m-1}'>\subset S_2( X_{\tau_{m-2}}\cap
{B_{k-2}},X_{\tau_{m-2}})$. Assume we found
vertices $a_m',a_{m-1}',...,a_l'$, $l>0$ spanning a simplex
in $S_2( X_{\tau_{l-1}}\cap
{B_{k-2}},X_{\tau_{l-1}})$, such that $\pi_k(a_i')=a_i$.
Then
we can find a vertex $a_{l-1}'\in
X_{\tau_{l-1}}\setminus {B_2( X_{\tau_{l-1}}\cap
{B_{k-2}},X_{\tau_{l-1}})}$ spanning together with $<
a_m',a_{m-1}',...,a_l'>$ a simplex in $X$.
Hence we can find
points $ a_m',a_{m-1}',...,a_0'\in S_k$ spanning a simplex in $X$
and satisfying $\pi_k(a_i')=a_i$. Then if
$z=\sum_{i=0}^m\lambda_ia_i$ for $\lambda_i>0$ such that
$\sum_{i=0}^m\lambda_i=1$, we have
$\pi_k(\sum_{i=0}^m\lambda_ia_i')=z$.
\end{proof}

\begin{lem}
\label{l7}
%\marginpar{[l7]}
Let the condition $R(w,X)$ hold
for every vertex $w$ of $X$. Then $\pi_k^{-1}(\rho)$ is connected for
every simplex $\rho$ of $S_{k-1}$ and for every $k\geq 2$.
\end{lem}
\begin{proof} 
If $\rho$ is a vertex then its preimage by the map
$\pi_k\colon S_k \to S_{k-1}$,
$\pi_k^{-1}(\rho)={\mr {span}}\lk
{\mr {vertices\; \; in}\;\; } X_{\rho}\setminus {B_2( X_{\rho}\cap
{B_{k-2}},X_{\rho})}\rik$  is nonempty, by Lemma \ref{l6} and it is
connected by $R(\tau,X)$.

By surjectivity of $\pi_k$ it is now enough to show the following. 
Let $\tau=<a_0,a_1,...,a_m>$ be a simplex of $S_{k-1}'$
such that $a_i$ is the barycenter of an $i$-simplex of
$S_{k-1}$. Then every point $p\in \pi_k^{-1}(\rho)$ can be
connected to $\pi_k^{-1}(a_0)$ by a path in
$\pi_k^{-1}(\tau)$.

The proof goes via induction on $m$.

For $m=0$ it is trivial.

Let $m>0$. W.l.o.g. we can assume that
$\pi_k (p)=\sum_{i=0}^{m}\lambda_ia_i$ with
$\lambda_m>0$ and $\sum_{i=0}^{m}\lambda_i=1$, $\lambda_i\geq 0$.
Let $I\subset \lk i | \lambda_i\neq 0\rik$ 
By definition of $\pi_k$ there exist
$(a_j')_{j\in I}$ such that $\pi_k(a_j')=a_j$,
$\sum_{j\in I}\lambda_j a_j'=p$.
Then the span of $(a_j')_{j\in I}$ is contained
in $\pi_k^{-1}(\tau)$ and hence we can connect
$p$ to $a_m'$ inside $\pi_k^{-1}(\tau)$.
Following {\it Case 2} in the proof of
Lemma \ref{l6} we can then find a vertex
$a_{m-1}''\in S_k$ such that
$\pi_k (a_{m-1}'')=a_{m-1}$ and $a_m'$
and $a_{m-1}''$ span an edge. Then
$<a_m',a_{m-1}''>\subset \pi_k^{-1}(\tau)$
and, by induction assumptions, we can connect
$p$ to $\pi_k^{-1}(a_0)$ by a path in
$\pi_k^{-1}(\tau)$.
\end{proof}

\begin{cor}
\label{c8}
%\marginpar{[c8]}
Let the condition $R(w,X)$ hold $X$
for every vertex $w$ of $X$. Then for every $k\geq 2$ and for any
connected subcomplex $K$ of $S_{k-1}$ its preimage $\pi_k^{-1}(K)$
is connected.
\end{cor}

\begin{tw}
\label{tw2}
%\marginpar{[tw2]}
Let $X$ be a finitely dimensional locally finite
$7$-systolic chamber complex such
that the condition $R(w,X)$ holds for every vertex $w$ of $X$.
Then the ideal boundary $\partial X$ of $X$ is connected.
\end{tw}
\begin{proof}
Observe that $S_1=X_v$ and thus it is
connected by $R(v,X)$. By Corollary \ref{c8} if $S_{k-1}$
is connected then
$S_k$ is connected too. Hence $\partial X$ as an inverse limit of
continua is a continuum.
\end{proof}

\begin{tw}
\label{tw3}
%\marginpar{[tw3]}
Let $X$ be a locally finite
$7$-systolic chamber complex of finite dimension $n\geq 3$.
Assume that
the
link $X_{\kappa}$ is connected, for every simplex $\kappa$ of $X$
of codimension greater then one, and $X_{\sigma}\setminus
\stackrel{\circ}{B_i}(\rho,X_{\sigma})$, $i=2,3$ is connected for every
codimension two simplex $\sigma$ of $X$ and every simplex $\rho$
of its link $X_{\sigma}$.
Then the ideal boundary $\partial X$ of $X$ is connected
and has no local cut points.
\end{tw}
\begin{proof}
Connectedness of the boundary follows from Lemma \ref{l0.5}
and Theorem \ref{tw2}.

Now we show there are no local cut points in $\partial X$.
If a point
$x\in \partial X$ disconnects an open connected set $U\subset
\partial X$ then it disconnects every open connected $V\subset U$.
Hence it disconnects every connected subset $W\subset U$ whose interior
contains $x$. Thus, to prove the Lemma, it is enough to show that for
a given point $x\in \partial X$ and its open neighbourhood $U$ there
exists a connected set $W$ with $x\in {\mr {int} \;}W\subset W\subset U$ such
that $W\setminus \lk x \rik$ is connected.

Let us define, for a natural number $k$, a map
$\pi^{\infty}_k\colon
\partial X \to S_k$ as a projection from the inverse limit
$\partial X$ to the element $S_k$ of the inverse system $\lk
S_i,\pi_i \rik$. By the definition of the topology on $\partial X$
we can find a natural number $k$ large enough so that if $\tau$
is a simplex of $S_k$ containing $\pi^{\infty}_k(x)$ then
$W=(\pi_k^{\infty})^{-1}({B_2(\tau,S_k)})\subset U$. We claim $W$
is as desired.

First observe that
$(\pi_k^{\infty})^{-1}(\stackrel{\circ}
{B_2}(\tau,S_k))\subset W$ is open and
contains $x$. Moreover ${B_2(\tau,S_k)}$ is a connected
subcomplex of $S_k$ and hence, by Corollary \ref{c8} the inverse
system $\lk W_l=\pi_l^{-1}(...
(\pi_{k+1}^{-1}
({B_2(\tau,S_k)}))...),\pi_l|_{W_l}
\rik _{l=k+1}^{\infty}$ consists of
continua and its inverse limit $W$ is a continuum.

Now we
show that every two points $y,z\in W\setminus\lk x \rik$ are
connected by a continuum within $W\setminus\lk x \rik$. Again by
the definition of the topology on $\partial X$ we can find $m$ big
enough such that there exists a vertex $w\in S_m$ such
that $\pi_m^{\infty}(x)\in  {\stackrel{\circ}
{B_1}(w,S_m)}$ and
$y,z\notin {\stackrel{\circ}{B_1}(w,S_m)}$.
Since $\pi_m^{\infty}(W)=\pi_{m}^{-1}(...(\pi_{k+2}^{-1}
(\pi_{k+1}^{-1}({B_2(\tau,S_k)})))...)$ is a
connected subcomplex of $S_m$
then, if $S_1(w,S_m)=(S_m)_w$ is
connected,
we can find a continuum $K\in
\pi_m^{\infty}(W)\setminus {\stackrel{\circ}{B_1}(w,S_m)}$
connecting $y$ and $z$.
Then $(\pi_m^{\infty})^{-1}(K)$ is a continuum (as an inverse
limit of continua) in $W$ missing $x$ and containing $y$ and $z$.

Thus to finish the proof we have to show that
$S_1(w,S_m)=(S_m)_w$ is connected.
Observe that for every simplex $\sigma$ of $X_w$ the link
of $\sigma$ in $X_w$ is the link of $\sigma \ast w$ in $X$.
Hence (compare Lemma \ref{l0.5} and its proof) the link
$X_w$ is a $7$-large chamber complex such that the condition
$R(z,X_w)$ holds for every vertex $z$, provided $X_w$ has dimension
above two. Let $\rho=X_w \cap S_{m-1}$.
Then, by Lemma \ref{l0.3}, we get $S_1(w,S_m)=(S_m)_w=
X_w\cap S_m=S_1(\rho, X_w)$. Since $\partial \rho$ is connected
and balls of small radii in $6$-large complexes are isomorphic
with the ones in its universal covers (cf. Lemma \ref{l0.15})
we get, by Corollary \ref{c8} that $S_1(\rho, X_w)$ is connected.
\end{proof}

\begin{cor}
\label{c9}
%\marginpar{[c9]}
Let $X$ be a locally finite normal $7$-systolic pseudomanifold
of finite dimension at least $3$. Then its ideal boundary
$\partial X$ is connected and has no local cut points.
\end{cor}
\begin{proof}
One-dimensional links in normal manifolds are circles.
\end{proof}

\begin{tw}
\label{tw4}
%\marginpar{[tw4]}
Let $G$ be a group acting geometrically by automorphisms on a
$7$-systolic chamber complex $X$ of dimension $n\geq 3$.
Assume that
the
link $X_{\kappa}$ is connected, for every simplex $\kappa$ of $X$
of codimension greater then one, and $X_{\sigma}\setminus
\stackrel{\circ}{B_i}(\rho,X_{\sigma})$,
$i=2,3$ is connected for every
codimension two simplex $\sigma$ of $X$ and every simplex $\rho$
of its link $X_{\sigma}$.
Then $G$ does not split essentially, as an amalgamated product or as an
$HNN$-extension, over a finite nor two-ended group.
\end{tw}
\begin{proof}
This follows from Theorem \ref{tw3}.

By Stallings theorem \cite{S}, $G$ does not split over a finite
group---compare also Gromov \cite[remarks after Proposition 3.2.A]{G},
Coornaert--Delzant--Papadopoulos \cite[Exercise 4) in Chapter 2]{CDP} and
Ghys and de la Harpe \cite[Proposition 17 in Chapter 7.5]{GdlH}.

By Bowditch \cite[Theorem 6.2]{B}, $G$ does not split essentially over
a two-ended group.
\end{proof}

\begin{cor}
\label{c10}
%\marginpar{[c10]}
A group acting geometrically by automorphisms on a
locally finite normal $7$-systolic pseudomanifold
of dimension at least $3$ does not split essentially,
as an amalgamated product or as an
$HNN$-extension, over a finite nor two-ended group.
\end{cor}

\begin{rem}
{\bf 1)} A systolic group acting on a $7$-systolic pseudomanifold of
dimension at least $3$ can split over a surface group (this remark
is due to J. {\' S}wi{\c a}tkowski). To see this take two
isomorphic closed $3$-dimensional $7$-large pseudomanifolds with
links of vertices being closed surfaces (such spaces exist by
Januszkiewicz and {\' S}wi{\c a}tkowski
\cite[Corollary 19.3 p. (1) and its proof]{JS1}).
Consider complement of an open residue of a given
vertex in each of them. The link of the vertex is a convex
subcomplex of the complement and hence the union of both complements
along that links is $7$-large. Thus the fundamental group of the
sum splits over the fundamental group of a link which is a
surface.

{\bf 2)} As noticed in \cite[Section 5]{O} most of the examples
of systolic groups---except automorphism groups of complexes
of dimension at most two---constructed in \cite{JS1} as fundamental
groups of some extra-tileable simplices of groups
satisfy assumtions of Theorems \ref{tw3} and \ref{tw4}.

At the moment the only $7$-systolic groups of virtual
cohomological dimension above two known to us are the groups
acting on normal $7$-systolic pseudomanifolds, constructed
in \cite{JS1}.
\end{rem}

{\bf Question.}

Can groups acting geometrically on normal
$7$-systolic pseudomanifolds of dimension at least $3$ split
over free non-abelian groups ?

%%%%%%%%%%%%%%%%%%%%   End of main body of article
%
%                             References
%
%   BiBTeX users uncomment the following line:
%
%\bibliographystyle{gtart}
%


\begin{thebibliography}
\bibitem[BeK]{BK} Benakli N., Kapovich I., Boundaries of hyperbolic
groups. Combinatorial and geometric group theory (New York,
2000/Hoboken, NJ, 2001), 39--93, Contemp. Math., 296, Amer. Math.
Soc., Providence, RI, 2002.
\bibitem[BesM]{BesM} Bestvina M., Mess G., {\it The boundary of
negatively curved groups},  J. Amer. Math. Soc.  4  (1991),  no. 3,
469--481.
%\bibitem[Bj]{Bj} A. Bj\"orner, {\it Nerves, fibers and homotopy groups},
%J. Combin. Theory Ser. A  102  (2003),  no. 1, 88--93.
%\bibitem[Bor]{Bo} K. Borsuk, {\it Theory of retracts},
%Monografie Matematyczne, Tom 44 Pa\'nstwowe Wydawnictwo Naukowe, Warsaw, 1967.
\bibitem[Bow]{B} Bowditch B. H., {\it Cut points and canonical
splittings of hyperbolic groups}, Acta Math. 180 (1998), no. 2,
145-186,
%\bibitem[B]{B} K. S. Brown, {\it Buildings}, Springer-Verlag, 1988,
\bibitem[BrH]{BH} Bridson M. R., Haefliger A., {\it Metric spaces of
non-positive curvature}, Grundlehren der Mathematischen
Wissenschaften , 319, Springer-Verlag, Berlin, 1999.
\bibitem[CDP]{CDP} Coornaert M., Delzant T., Papadopoulos A.,
{\it G{\' e}om{\' e}trie et th{\' e}orie des groupes},
Springer-Verlag, Berlin, 1990, Les groupes hyperboliques de
Gromov.
%\bibitem[Dav]{Dav} M. W. Davis, {\it Groups generated by reflectionsand
%aspherical manifolds not covered by Euclidean space}, Ann. of
%Math. ({\bf 2}) 117 (1983), 293-324.
\bibitem[Da]{Da} Daverman R. J., {\it Hereditarily aspherical
compacta and cell-like maps}, Topology Appl. 41 (1991),
no. 3, 247--254.
\bibitem[DaDr]{DD2} Daverman R. J., Dranishnikov A. N.,
{\it Cell-like maps and aspherical compacta},  Illinois J. Math.
40  (1996),  no. 1, 77--90.
%\bibitem[Dav]{Dav} Davis M. W., {\it The geometry and topology
%of Coxeter groups}, book in preparation.
\bibitem[DyO]{DO} Dymara J., Osajda D., {\it Boundaries of
right-angled hyperbolic buildings}, Fund. Math., to appear.
%\bibitem[G]{G} R. Geoghegan, {\it Topological methods in group
%theory}, book in preparation.
\bibitem[GdlH]{GdlH} Ghys {\' E}., de la Harpe P. (eds.),
{\it Sur les groupes hyperboliques d'apr{\` e}s Mikhael Gromov},
Birkh{\" a}user Boston Inc., Boston, MA, 1990.
\bibitem[G]{G} Gromov M., {\it Hyperbolic groups}, Essays in group
theory, Springer, New York, 1987, 75-263.
\bibitem[H]{H} Haglund F., {\it Complexes simpliciaux hyperboliques
de grande dimension}, preprint, Prepublication Orsay 71, 2003.
%\bibitem[HS]{HS} F. Haglund, J. \'Swi{\c a}tkowski,
\bibitem[Jak]{Jak} Jakobsche W., {\it Homogeneous cohomology
manifolds which are inverse limits}, Fund. Math. 137 (1991),
no. 2, 81--95.
\bibitem[JS1]{JS1} Januszkiewicz T., \'Swi{\c a}tkowski J.,
{\it Simplicial nonpositive curvature},  Publ. Math. Inst. Hautes 
�tudes Sci,  No. 104  (2006), 1--85.
\bibitem[JS2]{JS2} Januszkiewicz T., \'Swi{\c a}tkowski J., {\it
 Filling invariants of systolic complexes and groups},  Geom. Topol. 
11  (2007), 727--758.
\bibitem[JS3]{JS3} Januszkiewicz T,  \'Swi{\c a}tkowski J., {\it
Nonpositively curved developments of billiards}, submitted.
\bibitem[KaKl]{KaKl} Kapovich M., Kleiner B., {\it Hyperbolic groups
with low-dimensional boundary},  Ann. Sci. \'Ecole Norm. Sup. (4)
33  (2000),  no. 5, 647--669.
\bibitem[O]{O} Osajda D., {\it Connectedness at infinity of systolic complexes and
groups}, Groups, Geometry and Dynamics 1 (2007), 183--203.
\bibitem[PS]{PS} Przytycki P., \'Swi{\c a}tkowski J., {\it
Flag-no-square triangulations and Gromov boundaries in dimension
3}, submitted.
\bibitem[St]{S}
Stallings J. R., {\it On torsion-free groups with infinitely many ends},
Ann. of Math. (2)  88 (1968), 312--334.
\bibitem[W]{W} Wise D. T., {\it Sixtolic complexes and their fundamental groups},
in preparation.
\bibitem[Z]{Z} Zawi{\' s}lak P., {\it Pontryagin sphere as a boundary of 
systolic group}, in preparation.


\end{thebibliography}
\end{document}